\documentclass[11pt]{article}%
   
\usepackage{amsmath,enumerate}
\usepackage{amsfonts}
\usepackage{amssymb}
\usepackage{xcolor}

\newtheorem{theorem}{Theorem}[section]

\newtheorem{lemma}[theorem]{Lemma}

\newtheorem{proposition}[theorem]{Proposition}

\newtheorem{question}[theorem]{Question}

\newenvironment{proof}[1][Proof]{\noindent\textbf{#1.} }
{\hfill \ \rule{0.5em}{0.5em}}

\begin{document}

\title{Minimizing the number of complete bipartite graphs in a $K_s$-saturated graph}
\author{Beka Ergemlidze\thanks{Department of Mathematics and Statistics, University of South Florida,
Tampa, Florida, U.S.A. E-mail: \texttt{beka.ergemlidze@gmail.com}}\and Abhishek Methuku\thanks{School of Mathematics, University of Birmingham, United Kingdom. E-mail: \texttt{abhishekmethuku@gmail.com}. Research is supported by the EPSRC grant number EP/S00100X/1 (A.~Methuku).} \and Michael Tait\thanks{Department of Mathematics \& Statistics, Villanova University, U.S.A. E-mail: \texttt{michael.tait@villanova.edu}. Research is supported in part by National Science Foundation grant DMS-2011553.} \and Craig Timmons\thanks{Department of Mathematics and Statistics, California State University Sacramento, U.S.A. E-mail: \texttt{craig.timmons@csus.edu}. Research is supported in part by Simons Foundation Grant \#359419.}}
%
\maketitle
\begin{abstract}
A graph $G$ is $F$-saturated if it contains no copy of $F$ as a subgraph but the addition of any new edge to $G$ creates a copy of $F$. We prove that for $s \geq 3$ and $t \geq 2$, 
the minimum number of copies of $K_{1,t}$ in a 
$K_s$-saturated graph is $\Theta ( n^{t/2})$.  More precise results are obtained when $t = 2$ where the problem is related to 
Moore graphs with diameter 2 and girth 5.    
We prove that for $s \geq 4$ and $t \geq 3$, the minimum number of copies of $K_{2,t}$ in an $n$-vertex 
$K_s$-saturated graph is at least $\Omega( n^{t/5 + 8/5})$ and at most $O(n^{t/2 + 3/2})$.  
These results answer a question of 
Chakraborti and Loh.
General estimates on the number of copies of $K_{a,b}$ in a $K_s$-saturated graph are also obtained, but finding an asymptotic formula
remains open.  
\end{abstract}


\section{Introduction}

Let $F$ be a graph with at least one edge.
A graph $G$ is \emph{$F$-free} if $G$ does not 
contain $F$ as a subgraph.
The study of $F$-free graphs is central
to extremal combinatorics.  
Tur\'{a}n's Theorem, widely considered 
to be a cornerstone result in graph theory, 
determines the 
maximum number of edges in an $n$-vertex $K_s$-free
graph.  An interesting class of $F$-free graphs
are those that are maximal with respect to the 
addition of edges.  We say that a graph $G$ 
is \emph{$F$-saturated} if $G$ is $F$-free but the addition of an edge joining any pair of nonadjacent vertices of $G$ creates a copy of $F$.
The function $\textup{sat}(n , F)$ is the 
\emph{saturation number of $F$}, and is defined 
to be the minimum number of edges in an
$n$-vertex $F$-saturated graph.  In some 
sense, it is dual to the Tur\'{a}n function 
$\textup{ex}(n , F) $ which is the maximum 
number of edges in an $n$-vertex $F$-saturated graph.  

One of the first results on graph 
saturation is a theorem of 
Erd\H{o}s, Hajnal, and Moon \cite{ehm}
which determines the saturation 
number of $K_s$.  
They proved that for $2 \leq s \leq n$, 
there is a unique $n$-vertex $K_s$-saturated 
graph with the minimum number of edges.
This graph is the join of a complete 
graph on $s-2$ vertices and an independent 
set on $n - s + 2$ vertices,  
denoted $K_{s-2} + \overline{K_{n -s + 2}}$.  
The Erd\H{o}s-Hajnal-Moon Theorem 
was proved in the 1960s and since then, graph 
saturation has developed into its 
own area of extremal combinatorics.
We recommend the survey of Faudree, Faudree, and Schmitt \cite{sat survey}
as a reference for 
history and significant results in graph saturation.

The function $\textup{sat}(n, F)$ 
concerns the minimum number of edges in an $F$-saturated graph.  More generally, one can ask 
for the minimum number of copies of $H$ in an 
$n$-vertex $F$-saturated graph.  
Let us write $\textup{sat}(n , H , F)$ 
for this minimum.  This function was 
introduced in \cite{kmtt} and was 
motivated by the well-studied
generalized Tur\'{a}n function whose systematic study was initiated by Alon and Shikhelman~\cite{as}.  
Recalling that the Erd\H{o}s-Hajnal-Moon Theorem determines 
$\textup{sat}(n , K_s) = \textup{sat}(n , K_2 , K_s)$,
it is quite natural to study the generalized function 
$\textup{sat}(n , K_r , K_s)$, where 
$2 \leq r < s$.  Answering a question of Kritschgau, Methuku, Tait and Timmons \cite{kmtt},
Chakraborti and Loh 
\cite{cl} 
proved that for every $2 \leq r < s$, there is a constant $n_{r,s}$ such that 
for all $n \geq n_{r,s}$, 
\[
\textup{sat}(n , K_r , K_s ) = ( n - s + 2) \binom{s-2}{r-1} + \binom{s-2}{r}.
\]
Furthermore, they showed that 
$K_{s-2} + \overline{K_{ n -s + 2} }$ 
is the unique graph that minimizes the number of copies of $K_r$ among 
all $n$-vertex $K_s$-saturated graphs for $n \geq n_{r,s}$.  They proved a similar result 
for cycles where the critical point is that $K_{s-2} + \overline{K_{n - s + 2}}$ is
again the unique graph that minimizes the number of copies of $C_r$ among 
all $n$-vertex $K_s$-saturated graphs for $n \geq n_{r,s}$  under some assumptions on $r$ in relation to $s$ (see
Theorem \ref{c and loh} below). 
Chakraborti and Loh then asked the following question (Problem 10.5 in \cite{cl}).

\begin{question}
Is there a graph $H$ for which $K_{s-2} + \overline{K_{n - s + 2}}$ 
does not (uniquely) minimize the number of copies 
of $H$ among all $n$-vertex $K_s$-saturated graphs for 
all large enough $n$? 
\end{question}
 
Here we answer this question positively and show that 
there are graphs $H$ for which 
$K_{s-2} + \overline{K_{n - s + 2}}$ is \emph{not} the 
unique extremal graph.  

We begin by stating our first two results, 
Theorems \ref{k12 vs ks} and 
\ref{k1t vs ks}, where $H = K_{1,t}$.  Together, they
demonstrate a change in behaviour between the cases 
$H = K_{1,2}$ and $H = K_{1,t}$ with $t > 2$.  

\begin{theorem}\label{k12 vs ks}
\rm (i) For $n \geq 3$, 
\[
 \binom{n}{2}- \frac{ n^{3/2} }{2}
 \leq \textup{sat}(n , K_{1,2} , K_3) \leq 
 \binom{ n - 1}{2}  .
\] 

\smallskip
\noindent
(\rm ii) For $n \geq s \geq 4$,   
\[
\textup{sat}(n , K_{1,2} , K_s) 
= (s - 2) \binom{ n - 1}{2} + ( n - s + 2) \binom{s - 2}{2}.
\]
Furthermore, $K_{s-2} + \overline{K_{n - s  + 2}}$ is 
the unique $n$-vertex $K_s$-saturated with minimum number of copies of $K_{1,2}$.
\end{theorem}

\begin{theorem}\label{k1t vs ks}
For integers $n \geq s \geq 3$ and $t \geq 3$, 
\[
\textup{sat}(n , K_{1,t} , K_s ) = \Theta ( n^{t/2 + 1} ).
\]
\end{theorem}
A consequence of Theorem \ref{k1t vs ks} is that 
if $s ,t \geq 3$ and $n$ is large enough in terms of $t$, 
$K_{s-2} + \overline{K_{n - s + 2}}$ does not 
minimize the number of copies of $K_{1,t}$ among 
all $n$-vertex $K_s$-saturated graphs.  Indeed, 
$K_{s-2} + \overline{K_{n - s + 2}}$ has 
$\Theta (n^t )$ copies of $K_{1,t}$.  
Interestingly, the special case of determining 
$\textup{sat}(n , K_{1,2} , K_3)$ is related to the existence of Moore graphs.  
This is discussed further in the Concluding Remarks section, but whenever 
a Moore graph of diameter 2 and girth 5 exists, 
this graph will have fewer copies of $K_{1,2}$ than 
$K_1 + \overline{K_{n-1}}= K_{1 , n - 1} $.  Thus, any potential 
result that determines $\textup{sat}(n , K_{1,2} , K_3)$
exactly would have to take this into account.  

The graph used to 
prove the upper bound of 
Theorem \ref{k1t vs ks} is a $K_s$-saturated graph 
with maximum degree at most $c_s n^{1/2}$.  This graph 
was constructed by Alon, Erd\H{o}s, Holzman, and
Krivelevich \cite{aehk} and it is structurally very different 
from $K_{s-2} + \overline{K_{n - s + 2} }$.  
Using this graph one can prove 
a more general upper bound that applies to 
any connected bipartite graph.  This will be stated 
in Theorem \ref{general ub} below.

Next we turn our attention to counting copies of $K_{2,t}$ (for $t \geq 2$)
in $K_s$-saturated graphs.  The graph $K_1 + \overline{K_{n - 1}}$ is $K_3$-saturated and $K_{2,t}$-free.
Thus, $\textup{sat}(n , K_{2,t} , K_3) = 0$ for 
all $ t \geq 2$.  For $t = 2$ and $s \geq 4$ 
Chakraborti and Loh 
\cite{cl} proved that 
\begin{equation}\label{k22 ub of cl}
\textup{sat}(n , K_{2,2} , K_s ) 
= 
( 1 + o(1)) \binom{s-2}{2} \binom{n}{2}.
\end{equation}
Observe that the graph $K_{s - 2} + \overline{K_{n - s + 2}}$ 
has 
\[
\binom{s-2}{2} \binom{n - s + 2}{2} 
+ \binom{s-2}{3} (n - s + 2)
+ \binom{s-2}{4}
\]
copies of $K_{2,2}$ and this gives 
the upper bound in (\ref{k22 ub of cl}).  
Now the focus of \cite{cl} was on 
counting complete graphs and counting 
cycles, so here the above result is stated 
in terms of $K_{2,2}$ but of course
$K_{2,2} = C_4$.  
However, it is important and 
relevant to this work to mention the following theorem of Chakraborti and Loh which shows that $K_{s-2} + \overline{K_{ n - s + 2}}$ minimizes the number of copies of $C_r$ in certain cases.

\begin{theorem}[Chakraborti and Loh \cite{cl}]\label{c and loh}
Let $s \geq 4$ and $r \geq 7$ if $r$ odd, and 
$r \geq 4 \sqrt{s-2}$ if $r$ is even.  
There is an $n_{r,s}$ such that 
for all $n \geq n_{r,s}$, the graph 
$K_{s-2} + \overline{K_{n - s + 2}}$ minimizes 
the number of copies of $C_r$ over all $n$-vertex 
$K_s$-saturated graphs.  Moreover, when $r \leq 2s - 4$, 
this is the unique extremal graph.
\end{theorem}

It is conjectured in \cite{cl} that $K_{s-2} + \overline{K_{n - s + 2}}$ is the unique 
graph that minimizes the number of copies of $C_r$ among 
all $K_s$-saturated graphs.  Currently
it is only known that $K_{s-2} + 
\overline{K_{n - s + 2}}$ minimizes the number 
of copies of $K_r$ (Erd\H{o}s-Hajnal-Moon for $r = 2$ and \cite{cl} 
for $r > 2$), and minimizes the number of copies of $C_r$ under 
certain assumptions (stated in Theorem \ref{c and loh}).  Theorem \ref{k1t vs ks} shows 
$K_{s-2} + \overline{K_{n - s + 2}}$ \textit{does not} minimize
the number of copies of $K_{1,t}$.
We extend this to $K_{a,b}$ with $1 \leq a + 1 < b$ 
using the following theorem.

\begin{theorem}\label{general ub}
Let $F$ be a connected bipartite graph with parts of size $a$ and $b$ with
$1 \leq a + 1 < b$.  If $s \geq 3$ be an integer, then
\[
\textup{sat}(n , F , K_s ) = 
\left\{
\begin{array}{ll}
0 & \mbox{if $a > s - 2$}, \\
O ( n^{ \frac{1}{2}  (a + b + 1 ) } ) & \mbox{if $a \leq s - 2$}
\end{array}
\right.
\]
where the implicit constant can depend on $a$, $b$, and $s$.  
\end{theorem}

Theorem \ref{general ub} naturally suggests 
the following question: how many copies of $K_{2,t}$ must 
there be in a $K_s$-saturated graph?  In this direction we prove the following.

\begin{theorem}\label{main lower bound theorem}
Let $s \geq 4$ and $t \geq 3$ be integers.  
There is a positive constant $C$ such that 
\[
\textup{sat}(n , K_{2,t} , K_s ) \geq C n^{t/5 + 8/5}.
\]
\end{theorem}
By Theorem \ref{general ub}, 
$\textup{sat}(n, K_{2,t} , K_s) \leq O_{s,t} ( n^{t/2 + 3/2})$ for $s \geq 4$ and $t \geq 3$, so that there is a gap in the exponent in the upper and lower bounds.

Saturation problems with restrictions on the degrees 
have also been well-studied. Duffus and Hanson \cite{DH} investigated triangle-saturated graphs with minimum degree $2$ and $3$. 
Day \cite{Day} resolved a 20 year old conjecture of Bollob\'as \cite{Bollobas} which asked for a lower bound on the number of edges in $K_s$-saturated graphs with minimum degree $t$. Gould and Schmitt \cite{GS} studied $K_2^t$-saturated graphs (where $K_2^t$ is the complete $t$-partite graph with parts of size $2$)  with a given minimum degree.
Furthermore, $K_s$-saturated graphs with restrictions on the maximum degree were studied in \cite{aehk, fs, P}. Turning to generalized saturation numbers, as a step towards generalizing Day's 
Theorem, Curry et.\ al.\ \cite{CCDT} proved 
bounds on the number of triangles in a $K_s$-saturated graph with minimum degree $t$.    
Motivated by these results we prove a lower bound on the number of copies of $K_{a,b}$ in $K_s$-saturated graphs in terms of its minimum degree. 

\begin{theorem}\label{using hypergraph lemma min degree}
Let $s \geq 4$ and $2 \leq a  < b$ be integers with $a \leq s- 2$.  If $G$ is an $n$-vertex $K_s$-saturated graph with 
minimum degree $\delta (G)$, then 
$G$ contains at least 
\[
c \left(  \dfrac{ n - \delta (G) - 1 }{ \delta (G)^{a-1} } \right)^{b/2}
\]
copies of $K_{a,b}$ for some constant $c = c(s,a,b) > 0$.    
\end{theorem}

Theorem \ref{using hypergraph lemma min degree} shows that if $0 \leq \alpha < \frac{1}{a-1}$ and 
$\delta (G) \leq \kappa n^{ \alpha } $ for some $\kappa > 0$, then $G$ contains 
at least $c n^{ b/2( 1 - \alpha ( a - 1) ) }$ copies of $K_{a,b}$.  In particular, when 
$\delta(G)$ is a constant, we obtain $\Omega (n^{t/2})$ copies of $K_{2,t}$.
This improves the lower bound of Theorem \ref{main lower bound theorem}, but comes at the 
cost of a minimum degree assumption. 

In the next subsection we give the notation that 
will be used in our proofs.  Section \ref{k1t section} 
contains the proofs of Theorems \ref{k12 vs ks} and 
\ref{k1t vs ks}.  
Section 3 contains the proofs of 
Theorems \ref{general ub}, \ref{main lower bound theorem}, 
and \ref{using hypergraph lemma min degree}.

\subsection{Notation}

For graphs $F$ and $G$, we write 
$\mathcal{N} (F , G)$ for the number of copies of 
$F$ in $G$.  For a graph $G$ and $x,y \in V(G)$, 
write $N(x)$ for the neighborhood of $x$, and 
$N(x,y)$ for $N(x) \cap N(y)$. More generally, 
if $X \subseteq V(G)$ and $v \in V(G)$, then 
$N( v , X)$ is the set of vertices adjacent to 
all of the vertices in $\{v \} \cup X$, and 
$N(X)$ is the set of  vertices adjacent to all vertices in $X$.  
We write $d(v) = | N(v) |$, $d(X) = | N(X) |$, 
and $d(v,X) = | N(v,X) |$.  
The set $N[v] = N(v) \cup \{v \}$ is the closed neighborhood 
of $v$. For a graph $G$, let $e(G)$ denote the number of edges in $G$.

For a hypergraph $\mathcal{H}$, 
$d_{ \mathcal{H} } (v)$ is the number of edges in $\mathcal{H}$
containing $v$.  Similarly, 
$d_{ \mathcal{H} } (X)$ 
and 
$d_{ \mathcal{H} } (v,X)$ 
is the number of edges in $\mathcal{H}$ containing
$X$ and $\{v \} \cup X$, respectively.


\section{Bounds on $\textup{sat}(n , K_{1,t} , K_s )$}\label{k1t section}

\subsection{Proof of Theorem \ref{k12 vs ks}}

Since the graph $K_{s-2} + \overline{K_{n - s + 2} }$ is $K_s$-saturated, by counting the number of copies of $K_{1,2}$ in it, we have
\begin{equation}\label{new ub for k12 from ehm}
\textup{sat}(n , K_{1,2} , K_s)  \leq (s - 2) \binom{ n - 1}{2} + ( n - s + 2) \binom{s - 2}{2}.
\end{equation}
In particular, if $s = 3$ we have $\textup{sat}(n , K_{1,2} , K_3) \leq \binom{n-1}{2}$. We now prove a matching lower bound up to an error term of order $O ( n^{3/2})$.  
Let $G$ be an $n$-vertex $K_3$-saturated graph.  
If $e (G) \geq \frac{  \sqrt{  n - 1  } n }{2}$, then for $n \ge 3$,
\begin{eqnarray*}
\mathcal{N} ( K_{1,2} , G) & = & \sum_{v \in V(G) } \binom{d(v)}{2} 
\geq n \binom{ 2 e(G) / n }{2} = \frac{ 2e(G)^2 }{n } - e(G) \\ 
& \geq &   \binom{n}{2} - \frac{ n^{3/2} }{2}.
\end{eqnarray*}
Now assume that $e(G) <  \frac{  \sqrt{  n - 1  } n }{2}$.  If $x$ and $y$ are not adjacent, then since $G$ is $K_3$-saturated, $x$ and $y$ must be joined by a path of length 2.  Hence, 
\begin{equation*}
\mathcal{N} ( K_{1,2} , G) \geq  e( \overline{G} ) =  \binom{n}{2} - e (G) 
\geq  \binom{n}{2}- \frac{  n^{3/2} }{2}.
\end{equation*}
This completes the proof of (i) of Theorem  \ref{k12 vs ks}.  To prove (ii) of Theorem  \ref{k12 vs ks}, it suffices to show that for $ n \geq s \geq 4$,
$$\textup{sat}(n , K_{1,2} , K_s)  \ge (s - 2) \binom{ n - 1}{2} + ( n - s + 2) \binom{s - 2}{2},$$
since \eqref{new ub for k12 from ehm} holds.
Let $G$ be an $n$-vertex $K_s$-saturated graph
with $ n \geq s \geq 4$.  Kim, Kim, Kostochka and O \cite[Theorem 2.1]{kkko} proved that 
\begin{equation}
\label{eq:KKKO}
\sum_{v \in V(G) } (d(v)+1)(d(v)+2-s) \geq (s-2)n(n-s+1).    
\end{equation}
It is easy to check that 
\begin{equation}
\label{eq:rearrange}
\sum_{v \in V(G) } (d(v)+1)(d(v)+2-s)=\sum_{v \in V(G) } (d(v)-1)d(v)+(4-s)\sum_{v \in V(G)}d(v)+(2-s)n.    
\end{equation}
Therefore, combining \eqref{eq:KKKO} and \eqref{eq:rearrange}, we have
\begin{equation}\label{kim-kost-o-thm21-crl}
\sum_{v \in V(G) } (d(v)-1)d(v)\ge (s-2)n(n-s+1)+(s-4)2e(G)+(s-2)n.
\end{equation}
By the Erd\H{o}s-Hajnal-Moon Theorem
$$\textup{sat}(n ,K_s)=(s-2)(n-s+2)+\binom{ s - 2}{2},$$
and $K_{s-2} + \overline{K_{n - s + 2}}$ is the unique 
$n$-vertex $K_s$-saturated with $\textup{sat}(n , K_s)$ edges.  
Thus, 
$$2e(G)\ge 2(s-2)(n-s+2)+2\binom{ s - 2}{2}=(s-2)(2n-s+1).$$
Plugging this into \eqref{kim-kost-o-thm21-crl} we get that if $s\ge 4$,
$$\sum_{v \in V(G) } (d(v)-1)d(v)\ge (s-2)n(n-s+1)+(s-4)(s-2)(2n-s+1)+(s-2)n.$$
Dividing through by 2 and simplifying the right-hand side yields
$$\sum_{v \in V(G) } \binom{d(v)}{2}\ge (s - 2) \binom{ n - 1}{2} + ( n - s + 2) \binom{s - 2}{2},$$
where equality holds only if $G = K_{ s- 2} + \overline{K_{n - s  + 2}}$.  This completes
the proof of Theorem \ref{k12 vs ks}.

\subsection{Proof of Theorem \ref{k1t vs ks} }

Now we prove a lower bound on the number of copies of $K_{1,t}$
in a $K_s$-saturated graph that gives the correct order of magnitude for all $t \geq 3$.  

\begin{proposition}\label{new prop A}
Let $n \geq s \geq 3$ and $t \geq 3$ be integers. Then
\[
\textup{sat}(n , K_{1,t} , K_s )  \geq 
\left(   \frac{ \sqrt{s-2} }{t} \right)^t n^{t/2 + 1} + O_{s,t} (n^{t/2}).
\]
\end{proposition}
\begin{proof}
Let $G$ be an $n$-vertex $K_s$-saturated graph.  Kim, Kim, Kostochka and O \cite[Theorem 1.1]{kkko} proved 
that 
\begin{equation}\label{new prop A eq1}
\sum_{v \in V(G) } d(v)^2 \geq ( n - 1)^2 ( s - 2) + ( s -2)^2 ( n - s + 2)
\end{equation}
and that equality holds if and only if $G$ is $K_{s-2} + \overline{K_{n - s + 2}}$, except for 
in the case that $s = 3$ where equality holds if and and only if $G$ is $K_1 + \overline{K_{n - 1}}$ or a 
Moore graph.  
By the Power Means Inequality,
\begin{equation}\label{new prop A eq2}
\sum_{ v \in V(G) } d(v)^2 \leq n^{1 - 2/t} \left( \sum_{v \in V(G) } d(v)^t \right)^{2/t}.
\end{equation}
Combining (\ref{new prop A eq1}) and (\ref{new prop A eq2}) with the inequality 
$\sum_{v \in V(G) } d(v)^t \leq t^t \sum_{v \in V(G) } \binom{ d(v) }{t}$ and rearranging, we obtain that $\mathcal{N} ( K_{1,t} , G )$ is equal to
\[
 \sum_{ v \in V(G) } \binom{ d(v) }{t}  
\geq 
\frac{    (  ( n - 1)^2 ( s - 2) + ( s - 2)^2 ( n - s + 2) )^{t/2} }{ t^t n^{t/2 - 1} } = \left(   \frac{ \sqrt{s-2} }{t} \right)^t n^{t/2 + 1} + O_{s,t} (n^{t/2}).
\]
This completes the proof of Proposition \ref{new prop A}.  
 \end{proof}

\begin{proposition}\label{lb k1t vs ks}
Let $s \geq 3$ and $t \geq 3$ be integers.
For sufficiently large $n$, 
\[
\textup{sat}(n , K_{1,t} , K_s ) 
\leq \frac{ c_s^t n^{t/2 + 1} }{ t!}
\]
where $c_s$ is a constant depending only on $s$.
\end{proposition}
\begin{proof}
By a result of Alon, Erd\H{o}s, Holzman, and Krivelevich, for 
each $s \geq 3$ and sufficiently large $n$, there is a $K_s$-saturated graph $G$ 
with maximum degree $c_s \sqrt{n}$ (the constant $c_s$ satisfies
$c_s \rightarrow 2s$ as $s \rightarrow \infty$).  The number of copies of
$K_{1,t}$ in $G$ is then
\[
\sum_{v \in V(G) } \binom{d(v)}{t} \leq n \binom{ \Delta (G) }{t} 
\leq \frac{ c_s^t n^{t/2 + 1}}{t!}.
\]
\end{proof}

\begin{proof}[Proof of Theorem \ref{k1t vs ks}]
Theorem \ref{k1t vs ks} follows 
immediately from Propositions \ref{new prop A} and 
\ref{lb k1t vs ks}.
\end{proof}


\section{Bounds on $\textup{sat}(n , K_{2,t} , K_s )$ with $s \geq 4$ and $t \geq 3$}


\subsection{Upper bound on $\textup{sat}(n , K_{2,t} , K_s)$}

We begin this section with a basic lemma on counting 
copies of a graph $F$ in a graph $G$ with maximum degree $\Delta$.  
It is likely that this lemma, as well as Lemma \ref{Kab to F}, are known.

\begin{lemma}\label{general F lemma}
Let $F$ be a connected bipartite graph with parts of size $a$ and $b$.  If $G$ is an $n$-vertex 
graph with maximum degree $\Delta$, then 
\[
\mathcal{N}( F , G) \leq n \Delta^{a + b - 1}.
\]
\end{lemma}
\begin{proof}
We will prove the lemma by counting the number of possible embeddings of $F$ in $G$.
Let $d$ be the diameter of $F$, 
and $x$ be a vertex in $F$.  
For $0 \leq i \leq d$, let $N_i (x)$ be the set of vertices at distance $i$ from $x$ in $F$.  
We count embeddings of $F$ in $G$  
by starting with the vertex $x$, and then 
proceeding through $N_1(x)$, then $N_2(x)$ and so on.  
There are $n$ ways to choose a vertex in $G$ that corresponds to $x$.  Suppose that $v_x$ is the chosen vertex 
in $G$.  The vertices in $G$ 
corresponding to those in $N_1(x)$ must be neighbors of $v_x$ in $G$ and so there are at most 
$\Delta^{ |N_1(x)| }$ possibilities.  This process is then repeated on $N_2(x)$, $N_3(x)$, and so on.  
The crucial point is that each time a vertex of $F$ is embedded in $G$, it is a neighbor (in $G$) of a
previously embedded vertex (from $F$).  Therefore, the number of possible embeddings of $F$ in $G$ 
is at most 
\[
n \Delta^{ |N_1(x)| } \Delta^{ |N_2(x) | } \cdots \Delta^{ | N_d(x) | } = n \Delta^{a + b - 1}.
\]
Here we have used the assumption that since $F$ is a connected graph with diameter $d$, we have 
the partition
\[
\{ x \} \cup N_1 (x) \cup N_2(x) \cup \dots \cup N_d(x) = V(F).
\]
\end{proof}

\begin{lemma}\label{Kab to F}
Let $F$ be a connected bipartite graph with parts of size $a$ and $b$.  For any $n$-vertex graph $G$,
\[
\mathcal{N} ( K_{a,b} , G ) \leq \mathcal{N} ( F , G ).
\]
\end{lemma}
\begin{proof}
If $G$ has no $K_{a,b}$, then the lemma is trivial.  Suppose $K$ is a copy of $K_{a,b}$ in $G$.  
Then, since $F$ is a subgraph of $K_{a,b}$, we have that $F$ is a subgraph of $K$ so $G$ has a copy of $F$.
Moreover, since any two different copies of $K_{a,b}$ have different vertex sets, they give rise to different copies of $F$.  Thus, for each copy of
$K_{a,b}$ in $G$ we obtain a copy of $F$, and no copy of $F$ will be obtained twice in this way. This proves Lemma \ref{Kab to F}.
\end{proof}

\medskip

We are now ready to prove Theorem \ref{general ub}.  

\begin{proof}[Proof of Theorem \ref{general ub}]
If $a > s-2$, then $K_{s-2} + \overline{K_{n - s + 2}}$ is $K_s$-saturated with no copies of $F$.  Indeed, a copy of $F$ would need 
at least $a$ vertices from the $K_{s-2}$, but $a > s - 2$. 

Now assume $a \leq s - 2$.  Let $G_q^s$ be the $K_s$-saturated graph constructed in \cite{aehk} where
$n$ (and thus $q$) is chosen large enough so that $b < \frac{q+1}{2}$.   There is a constant 
$c_s > 0$ such that $\Delta (G_q^s ) \leq c_s \sqrt{n}$.  By Lemma \ref{general F lemma},
the number of copies of $F$ in $G_q^s$ is at most $n c_s^{a+b-1} n^{ (1/2) ( a + b - 1)  } = c_s^{a+b+1} n^{ (1/2)(a + b + 1)}$.
 \end{proof}

\medskip

We conclude this subsection 
by showing that the graph $G_q^s$ used in the proof of Theorem \ref{general ub} cannot 
be used to further improve upon the upper bound of 
$O( n^{ \frac{1}{2} ( a + b + 1) })$
when $F = K_{a,b}$.  Since we 
are showing that $G_q^s$ cannot be used to improve the upper bound, we will be brief in our argument.  
We will use the same terminology as in \cite{aehk}, but  
one point at which we differ is the 
notation we use for a vertex.  A vertex 
in $G_q^s$ is determined by its level, place, type, and 
copy.  A vertex at level $i$, place $j$, type $t$, and 
copy $c$ will be written as 
\[
( ( i - 1) q + j , t , c).
\]
  
First, take $n$ large enough so that $b < \frac{q+1}{2}$.  
Choose a sequence $i_1, i_2 , \dots , i_a$ of 
levels 
with $1 \leq i_1 < i_2 < \dots < i_a \leq \frac{q+1}{2}$.  
Likewise, choose a sequence of $b$ levels 
$\frac{q+1}{2} \leq i_{a+1} < i_{a+2} < \dots < i_{a + b} \leq q+1$.
This can be done in 
$\binom{ \frac{q+1}{2} }{a} \binom{  \frac{q+1}{2} }{b}$ ways.  Next, choose a place $j_1 \in [q]$ which can 
be done in $q$ ways, and a type $t_1 \in [s-2]$ which can be done in $s-2$ ways.  Finally, choose a sequence of 
copies $1 \leq c_1 , c_2 , \dots , c_{a+b} \leq s-1$
arbitrarily.  This can be done in $(s-1)^{a+b}$ ways.
Using the definition of $G_q^s$, one finds that the $a$ vertices in the set
\[
\{  (  (i_z - 1)q + j_1 , t_1 , c_z ) : 1 \leq z \leq a \}
\]
are all adjacent to the $b$ vertices in the set
\[
\{  (  (i_z - 1)q + (j_1 + 1)_q  , 
(t_1 + 1)_{s-2} , c_z ) 
: a +1 \leq z \leq a + b  \}
\]
(here $(j_1+1)_q$ is the unique integer $\zeta$
in $\{1,2, \dots ,q \}$ for 
which $j_1 + 1 \equiv \zeta ( \textup{mod}~q)$, 
and $(t_1 + 1)_{s-2}$ is the unique integer $\zeta'$ in 
$\{1,2, \dots , s- 2 \}$ for which 
$t_1 + 1 \equiv \zeta ' ( \textup{mod}~s - 2)$).  
This gives a $K_{a,b}$ in $G_q^s$ and so 
the number of $K_{a,b}$ in $G_q^s$ is at least 
\[
\binom{ (1/2)( q  + 1) }{a} \binom{ (1/2) ( q + 1) }{b} 
q ( s - 2) ( s - 1)^{a+b} 
\geq 
C_{s,a,b} q^{ a + b + 1 } 
\geq 
C n^{ (1/2) ( a + b + 1) }.
\]
By Lemmas \ref{Kab to F} and \ref{general F lemma}, $G_q^s$ is a $K_s$-saturated $n$-vertex 
graph with 
$\Theta_{s,a,b} ( n^{ (1/2) (a + b +1 ) })$
copies of $K_{a,b}$.


\subsection{Lower bound on $\textup{sat}(n , K_{2,t} , K_s )$}

First we prove Theorem \ref{main lower bound theorem}.

\begin{proof}[Proof of Theorem \ref{main lower bound theorem}]
Let $G$ be a $K_s$-saturated graph on $n$ vertices. Note that we can assume
\begin{equation}
\label{eq:lowerboundeG}    
e(G)\le \frac{n^{\frac74}}{10}.  
\end{equation}
Otherwise, 
a theorem of Erd\H{o}s and Simonovits \cite{erdos-simonovits} implies that
there is a positive constant $\gamma$ such that  
\begin{equation}
\label{eq:supersat}
\mathcal N(K_{2,t},G)\ge\gamma\frac{e(G)^{2t}}{n^{3t-2}} = \Omega(n^{\frac{t}{2}+2}),
\end{equation}
proving Theorem \ref{main lower bound theorem}.

Let $K_4^-$ be the graph consisting of $4$ vertices and $5$ edges obtained by removing an edge from $K_4$. For a copy of $K_4^-$ with vertices $x,y,u,v$, where $uv\notin E(G)$, let $xy$ be called the \emph{base edge} of this $K_4^-$.
We estimate the number of copies of $K_4^-$ in a $K_s$-saturated graph $G$.

For every $u,v$ with $uv\notin E(G)$ there is a set $S$ such that $S \subseteq N(u,v)$ and
$S$ induces a $K_{s-2}$ in $G$. Therefore, there are at least $\binom{s-2}{2}$ pairs $x,y\in S$ such that $u,v,x,y$ form a copy of $K_4^-$.
On the other hand, every $xy\in E(G)$ is the base edge of at most $\binom{d(x,y)}{2}$ copies of $K_4^-$ in $G$.

Therefore,
$$\sum_{xy\in E(G)}{\binom{d(x,y)}{2}}\ge  \mathcal{N} ( K_4^- , G) \ge \sum_{uv\in E(\overline{G})}\binom{s-2}{2}\ge e(\overline G)\overset{\eqref{eq:lowerboundeG}}{\ge} \frac{n^2}{4}.$$

Thus, there is a constant $c_t = c(t)$ such that the following holds:

$$\mathcal N(K_{2,t},G)\ge\sum_{xy\in E(G)}\binom{d(x,y)}{t}\ge \frac{1}{t^t} \sum_{xy\in E(G)}\left({\binom{d(x,y)}{2}}^{\frac t2}-t^t\right)\ge$$

$$ \frac{e(G)}{t^{t}}\left(\frac{\sum_{xy\in E(G)}{\binom{d(x,y)}{2}}}{e(G)}\right)^{\frac t2}-e(G)\ge \frac{(n^2/4)^{\frac t2}}{t^te(G)^{\frac t2-1}}-e(G)=\frac{(n^2/4)^{\frac t2}-t^te(G)^{t/2}}{t^te(G)^{\frac t2-1}}\overset{\eqref{eq:lowerboundeG}}{\ge}\frac{c_tn^t}{e(G)^{\frac t2-1}}.$$

Combining this with \eqref{eq:supersat} we get

$$\mathcal N(K_{2,t},G)\ge \min \{\gamma\frac{e(G)^{2t}}{n^{3t-2}} , \frac{c_tn^t}{e(G)^{\frac t2-1}}
\}.
$$
Let $e(G)=n^{\alpha}$, then 
$$\mathcal N(K_{2,t},G)\ge \min \left \{\gamma{n^{2\alpha t-3t+2}} ,c_tn^{t-\alpha t/2+\alpha}\right
\}.
$$
Choosing $\alpha=\frac{8t-4}{5t-2}$ and $C = \min\{\gamma, c_t\}$, we get the desired lower bound $C n^{\frac t5-\frac{16}{125t-10}+\frac{41}{25}}>C n^{\frac t5+\frac 85}$.
\end{proof}

Next we turn to the proof of Theorem \ref{using hypergraph lemma min degree}. We need the following lemma.

\begin{lemma}\label{general hypergraph lemma}
Let $s \geq 4$ and $2 \leq a \leq b$ be integers with $a \leq s - 2$.  Suppose 
that $G$ is an $n$-vertex $K_s$-saturated graph with vertex set $V$.  There is a constant $c = c ( s , a , b)$ such that for any $v \in V$, there are at least 
\[
c \left( \frac{ n- d(v) - 1 }{ d(v)^{a - 1} } \right)^{ b / 2}
\]
copies of $K_{a,b}$ containing $v$.  
\end{lemma}
\begin{proof}
Let $v \in V $.  For each $u \in V \backslash N[v]$, there is a set $S_u \subset N(v)$ such that 
$S_u$ induces a $K_{s-2}$ in $G$.  Fix such an $S_u$ and define an 
$(s - 1)$-uniform hypergraph $\mathcal{H}$ to have vertex set $V \backslash \{ v\}$, 
and edge set $E ( \mathcal{H} ) = \{ \{u \} \cup S_u : u \in V \backslash N[v] \}$.  
By construction, $\mathcal{H}$ has $n - d(v) - 1$ edges, each of which contains 
exactly one vertex from $V \backslash N[v]$ and $s - 2$ vertices from $N(v)$.  Also, no 
two edges of $\mathcal{H}$ contain the same vertex from $V \backslash N[v]$.  
In what follows, we will add the subscript $\mathcal{H}$ if we are referring to degrees in 
$\mathcal{H}$, and no subscript will be included if we are referring to degrees or neighborhoods in $G$.  

By averaging, there is a set $X \in \binom{ N(v) }{ a - 1}$ such that 
\[
d_{ \mathcal{H} } (X) \geq \frac{  \binom{s-2}{a-1}  ( n-  d(v) - 1 ) } {  \binom{ d(v) }{ a - 1} }.
\]
We then have 
\begin{equation}\label{general hypergraph lemma eq 1}
\sum_{ y \in N (v,X) } d_{ \mathcal{H} } ( y , X) 
\geq \frac{  d_{ \mathcal{H} } (X) }{ (s - 2) - |X| } \geq c_1 \frac{ n - d(v) - 1}{ d(v)^{a-1} }
\end{equation}
for some constant $c_1 = c_1 ( s , a) > 0$.  The number of $K_{a,b}$ with $X \cup \{y \}$ forming 
the part of size $a$ ($y$ is an arbitrary vertex from $ N(v,X) $) and $v$ in the part of size $b$ is at least 
\[
\sum_{ y \in N (v , X) } \binom{  d_{ \mathcal{H} } ( y , X) }{ b - 1} 
\geq d(v,X) \binom{   \frac{ c_1 ( n - d(v) - 1 ) }{ d(v,X) d(v)^{a-1} }  }{b - 1}
\geq
 \frac{ c_2 ( n - d(v) - 1)^{ b - 1} }{ d(v,X)^{b - 2}  d(v)^{ ( a- 1)( b - 1) } }.
 \]
Here we have used convexity, (\ref{general hypergraph lemma eq 1}), and $c_2 = c_2 ( s , a, b)$ is some 
positive constant.  

Recalling that $|X| = a -1$, there are $\binom{ d(v,X) }{ b}$ copies of 
$K_{a,b}$ where $\{v \} \cup X$ is the part of size $a$ and the 
part of size $b$ is contained in $N(v) \backslash X$.  Thus, for some constant $c_3 = c_3 ( s, a, b) > 0$, the number 
of $K_{a,b}$ that contain $v$ is at least
\[
\frac{c_3 ( n - d(v) - 1)^{ b - 1} }{ d(v,X)^{b - 2} d(v)^{ (a - 1)( b - 1) } } + c_3 d(v,X)^b.
\]
By considering cases as to which is this the bigger term in 
this sum, we find that in both cases, there are at least 
\[
c_3 \left(  \frac{ n - d(v) - 1}{  d(v)^{ a- 1}} \right)^{b /2} 
\]
copies of $K_{a,b}$ containing $v$.   
\end{proof}

\medskip

Applying Lemma \ref{general hypergraph lemma} to a vertex $v$ with $d(v) = \delta (G)$ proves 
Theorem \ref{using hypergraph lemma min degree}.


\section{Concluding Remarks}

An interesting open problem is determining
the minimum number of copies of $K_{1,2}$ in a $K_3$-saturated graph.  
There is a connection between this problem and 
Moore graphs with diameter 2 and girth 5.  It is easy to check that an $n$-vertex 
Moore graph with diameter 2 and girth 5 is $K_3$-saturated, and it is regular with degree $d = \sqrt{n-1}$ \cite{singleton}
so it contains $n \binom{d}{2} = n \binom{ \sqrt{n-1} }{2}$ copies of $K_{1,2}$, and for all $n \ge 3$, this value is less than $\binom{n-1}{2}$ which 
is the number of copies of $K_{1,2}$ in $K_1 + \overline{K_{n-1}} = K_{1,n-1}$.
Furthermore, one can duplicate vertices of 
a Moore graph and preserve the $K_3$-saturated property (where each duplicated vertex has the same neighborhood as the original vertex).
Duplicating a vertex of the Petersen graph will lead to 
an $11$-vertex $K_3$-saturated graph with 42 copies of $K_{1,2}$, 
but $K_{1,10}$ has 45 copies of $K_{1,2}$.  Starting 
from the Hoffman-Singleton graph, one can duplicate
a vertex up to 4 times and we can still have fewer
copies of $K_{1,2}$ compared to the number of copies of $K_{1,2}$ in $K_{1,n-1}$.  
Duplicating a single vertex is not necessarily the optimal way to minimize the number of copies of $K_{1,2}$, but the point 
is that there are other graphs besides the Moore graphs that have fewer copies of $K_{1,2}$ than the number of copies of $K_{1,2}$ in $K_{1,n-1}$. 

It would also be interesting to determine the 
order of magnitude of $\textup{sat}(n , K_{2,t} , K_s )$.  There 
is a gap in the exponents (which is discussed in the introduction)
and it would be nice to close this gap.  It is not clear if our lower or upper bound is closer to the 
correct answer.  

Another potential approach to studying $\textup{sat}(n, H, F)$ is via the random $F$-free process. This random process orders the pairs of vertices uniformly and then adds them one by one subject to the condition that adding an edge does not create a copy of $F$. The resulting graph is then $F$-saturated. This process was first considered in \cite{br, esw, rw, s} and has since been studied extensively.  If $X_{H, F}$ is the random variable that counts the number of copies of $H$ in the output of this process, then we have that $\textup{sat}(n, H, F) \leq \mathbb{E}(X_{H,F})$. It would be interesting to determine for which graphs $H$ and $F$ that this approach gives better bounds than the explicit constructions that are currently known.



\end{document}